\documentclass[oneside,a4paper]{article}
\usepackage{en-note}
\usepackage{indentfirst}
\hypersetup{
pdftitle={Alexandrov's embedding theorem},
pdfauthor={Nina Lebedeva and Anton Petrunin}
}

\begin{document}

\title{Alexandrov's embedding theorem}
\author{Nina Lebedeva and Anton Petrunin}
\date{}
\maketitle

\section{Introduction}

Intrinsic distance between two points on the surface of a convex polyhedron is defined as the length of a shortest curve on the surface between these points.

Recall that the sum of angles at the tip of a convex polyhedral angle is less than $2\cdot\pi$;
this statement can be found in a school textbook \cite[§~48]{kiselyov-3D}.

It is easy to see that the surface of a convex polyhedron is homeomorphic to the sphere.
Therefore the statements above imply that the surface of a convex polyhedron equipped with its intrinsic metric is an example of a \emph{polyhedral metric on the sphere with the sum of angles around each vertex at most $2\cdot\pi$};
a metric is called \emph{polyhedral} if the sphere admits a triangulation such that every triangle is congruent to a plane triangle.

Alexandrov's theorem states that the converse holds if one includes in the consideration \emph{twice covered polygons}.
In other words, we assume that a polyhedron can degenerate to a plane polygon;
in this case, its surface is defined as two copies of the polygon glued along their boundary.

Further, we assume that a polyhedron can degenerate to a plane polygon.

\begin{thm}{Alexandrov's theorem}
\begin{enumerate}[I.]
\item\label{thm:exist}
A polyhedral metric on the sphere is isometric to the surface of a convex polyhedron if and only if the sum of angles around each of its vertex is not greater than $2\cdot\pi$.

\item\label{thm:unique} 
Moreover, a convex polyhedron is defined up to congruence by the intrinsic metric on its surface.
\end{enumerate}

\end{thm}

A. D. Alexandrov has many remarkable theorems, but in our opinion, this theorem is the most remarkable.
At the same time, its proof is elementary;
it could be explained to anyone familiar with basic topology.

This theorem has many applications.
In particular, it is used in the proof of its generalization \cite{alexandrov-1948} that gives a complete description of intrinsic metrics on the sphere that are isometric to convex surfaces in the Euclidean space.
The latter statement is fundamental in a branch of modern mathematics --- the so-called \emph{Alexandrov geometry}.

The first part is central; it is called the \emph{existence theorem}.
The second part is called the \emph{uniqueness theorem}; it is a slight variation of Cauchy's theorem about polyhedrons.
(There is another uniqueness theorem of Alexandrov that generalizes Minkowski's theorem about  polyhedrons.)

According to the theorem, a convex polyhedron is completely defined by the intrinsic metric of its surface.
In particular, knowing the metric we could find the position of the edges.
However, in practice, it is not easy to do.
For example, the surface glued from a rectangle as shown on the diagram defines a tetrahedron.
Some of the glued lines appear inside facets of the tetrahedron and some edges (dashed lines) do not follow the sides of the rectangle.

{

\begin{wrapfigure}{r}{30mm}
\vskip-3mm
\centering
\includegraphics{mppics/pic-10}
\vskip-0mm
\end{wrapfigure}

The theorem was proved by A. D. Alexandrov in 1941 \cite{alexandrov-1941};
we will present a sketch of his proof.
A complete proof is nicely written by A. D. Alexandrov in his book~\cite{alexandrov}.
Yet another proof was found by Yu.~A.~Volkov in his thesis \cite{volkov};
it uses a deformation of three-dimensional polyhedral space.

}

\section{Space of polyhedrons and metrics}

\paragraph{Space of polyhedrons.}
Let us denote by $\Phi$ the space of all convex polyhedrons in the Euclidean space,
including polyhedrons that degenerate to a plane polygon.
Polyhedra in $\Phi$ will be considered up to a motion of the space, 
and the whole space $\Phi$ will be considered with the natural topology (an intuitive meaning of closeness of two polyhedrons should be sufficient).  

Further, denote by $\Phi_n$ the polyhedrons in $\Phi$ with exactly $n$ vertices.
Since any polyhedron has at least 3 vertices, the space $\Phi$ admits a subdivision into a countable number of subsets $\Phi_3,\Phi_4,\dots$

\paragraph{Space of polyhedral metrics.}
The space of polyhedral metrics on the sphere with the sum of angles around each point at most $2\cdot\pi$ will be denoted by $\Psi$.
The metrics in $\Psi$ will be considered up to an isometry, and the whole space $\Psi$ will be equipped with the natural topology (again, an intuitive meaning of closeness of two metrics is sufficient).

A point on the sphere with the sum of angles strictly less than $2\cdot\pi$ will be called an \emph{essential vertex}.
The subset of $\Psi$ of all metrics with exactly $n$ essential vertices will be denoted by $\Psi_n$.
It is easy to see that any metric in $\Psi$ has at least 3 essential vertices.
Therefore $\Psi$ is subdivided into countably many subsets
 $\Psi_3,\Psi_4,\dots$

\paragraph{From a polyhedron to its surface.}

Recall that the surface of a convex polyhedron is a sphere with a polyhedral metric such that the sum of angles around each point is at most $2\cdot\pi$.
Therefore passing from a polyhedron to its surface defines a map
\[\iota\:\Phi\to \Psi.\]

Note that the number of vertices of a polyhedron is equal to the number of essential vertices of its surface.
In other words, $\iota(\Phi_n)\subset \Psi_n$ for any $n\ge 3$.

\section{About the proof}

Using the notation introduced in the previous section, we can give the following more exact formulation of Alexandrov's theorem: 

\begin{thm}{Reformulation}
For any integer $n\ge 3$,
the map $\iota$ is a bijection from $\Phi_n$ to~$\Psi_n$.
\end{thm}

We sketch the original proof of A. D. Alexandrov.
It is based on the  construction of a one-parameter family of polyhedrons that starts at arbitrary polyhedron
and ends at an arbitrary polyhedron with its surface isometric to the given one.
This type of argument is called the \emph{continuity method}; it is often used in the theory of differential equations.

\medskip

The two parts of the first formulation will be proved separately.

\parit{Part \ref{thm:unique}.} Let us show that the map $\iota\:\Phi_n\to\Psi_n$ is injective;
in other words, a convex polyhedron is defined by the intrinsic metric on its surface up to a motion of the space.

The last statement is analogous to the Cauchy theorem about polyhedrons,
and the proof goes along the same lines. 

The Cauchy theorem states that facets of a polyhedron together with the gluing rule completely describe a convex polyhedron;
its proof is given in many classical popular texts \cite{aigner-zigler,dolbilin,tabacnikov-fuks}.

\medskip

\parit{Part \ref{thm:exist}.}
Let us prove that $\iota\:\Phi_n\to\Psi_n$ is surjective.
This part of the proof is subdivided into the following lemmas:

\begin{thm}{Lemma}
For any integer $n\ge 3$, the space $\Psi_n$ is connected.
\end{thm}

The proof of this lemma is not complicated, but it requires ingenuity;
it can be done by the direct construction of a one-parameter family of metrics in $\Psi_n$ that connects two given metrics.
Such a family can be obtained by а sequential application of the following construction and its inverse.

Let $M$ be a sphere with metric from $\Psi_n$.
Suppose $v$ and $w$ are essential vertices in $M$.
Let us cut $M$ along a shortest line from $v$ to $w$.
Note that the shortest line cannot pass thru an essential vertex of $M$.
Further, note that there is a three-parameter family of patches that can be used to patch the cut so that the obtained metric remains in $\Psi_n$;
in particular, the obtained metric has exactly $n$ essential vertices (after the patching, the vertices $v$ and $w$ may become inessential).

\begin{thm}{Lemma}
The map $\iota\:\Phi_n\to\Psi_n$ is open, 
that is, it maps any open set in $\Phi_n$ to an open set in $\Psi_n$.

In particular, for any $n\ge 3$, the image $\iota(\Phi_n)$ is open in~$\Psi_n$.
\end{thm}

This statement is very close to the so-called \emph{invariance of domain theorem};
the latter states that a continuous injective map between manifolds of the same dimension is open.

According to part \ref{thm:unique}, $\iota$ is injective.
The proof of the invariance of domain theorem can be adapted to our case since both spaces $\Phi_n$ and $\Psi_n$ are $(3\cdot n-6)$-dimensional and both look like manifolds, altho, formally speaking, they are \emph{not} manifolds.
In a more technical language, $\Phi_n$ and $\Psi_n$ have the natural structure of $(3\cdot n-6)$-dimensional \emph{orbifolds},
and the map $\iota$ respects the \emph{orbifold structure}.

We will only show that both spaces $\Phi_n$ and $\Psi_n$ are $(3\cdot n-6)$-dimensional.

Choose a polyhedron $P$ in $\Phi_n$.
Note that $P$ is uniquely determined by the $3\cdot n$ coordinates of its $n$ vertices.
We can assume that the first vertex is the origin, the second has two vanishing coordinates and the third has one vanishing coordinate; therefore, all polyhedrons in $\Phi_n$ that lie sufficiently close to $P$ can be described by $3\cdot n-6$ parameters.
If $P$ has no symmetries then this description can be made one-to-one;
in this case, a neighborhood of $P$ in $\Phi_n$ is a $(3\cdot n-6)$-dimensional manifold.
If $P$ has a nontrivial symmetry group, then this description is not one-to-one but it does not have an impact on the dimension of $\Phi_n$.

The case of polyhedral metrics is analogous.
We need to construct a subdivision of the sphere into plane triangles using only essential vertices.
By Euler's formula, there are exactly $3\cdot n-6$ edges in this subdivision.
Note that the lengths of edges completely describe the metric, and slight changes of these lengths produce a metric with the same property.

\begin{thm}{Lemma}
The map $\iota\:\Phi_n\to\Psi_n$ is closed;
that is, the image of a closed set in $\Phi_n$ is closed in $\Psi_n$.

In particular, for any $n\ge 3$, the set $\iota(\Phi_n)$ is closed in~$\Psi_n$.
\end{thm}

Choose a closed set $Z$ in $\Phi_n$.
Denote by $\bar Z$ the closure of $Z$ in $\Phi$; note that $Z=\Phi_n\cap \bar Z$.
Assume $P_1,P_2,\dots\in Z$ is a sequence of polyhedrons that converges to a polyhedron $P_\infty\in\bar Z$.
Note that $\iota(P_n)$ converges to $\iota(P_\infty)$ in $\Psi$.
In particular, $\iota(\bar Z)$ is closed in $\Psi$.

Since $\iota(\Phi_n)\subset \Psi_n$ for any $n\ge 3$, we have  $\iota (Z)=\iota(\bar Z)\cap \Psi_n$;
that is, $\iota (Z)$ is closed in $\Psi_n$. 

\medskip

Summarizing, $\iota(\Phi_n)$ is a nonempty closed and open set in $\Psi_n$, and $\Psi_n$ is connected for any $n\ge 3$.
Therefore, $\iota(\Phi_n)=\Psi_n$; that is, $\iota\:\Phi_n\z\to\Psi_n$ is surjective.
\qeds

\parbf{Acknowledgments.} We want to thank S. Alexander, Y. Burago, and J. Tsukahara for help. 
The authors were partially supported by RFBR grant 20-01-00070 and NSF grant DMS-2005279.

\sloppy
\printbibliography[heading=bibintoc]
\fussy

\end{document}


\title{Теорема Александрова о\\ вложении многогранников}
\author{Нина Лебедева и Антон Петрунин}
\date{}
\maketitle

\section{Введение}

Нас будет интересовать \emph{внутренняя метрика} (то есть метрика, задающаяся длинами кратчайших путей между точками) на поверхности выпуклого многогранника.

Напомним, что сумма углов при вершине многогранного выпуклого угла меньше $2\cdot \pi$; 
это утверждение можно найти в школьном учебнике \cite[§~325]{kiselyov}.

Нетрудно видеть, что поверхность выпуклого многогранника гомеоморфна сфере.
Из вышесказанного получаем, что поверхность выпуклого многогранника, наделённая естественной внутренней метрикой,
является примером \emph{многогранной метрики на двумерной сфере с суммой углов вокруг каждой вершины, не превосходящей $2\cdot\pi$}; 
метрика называется \emph{многогранной} если сфера допускает триангуляцию такую, что каждый её треугольник равен плоскому треугольнику.

Теорема Александрова гласит, что обратное верно, если включить в рассмотрение \emph{дважды покрытые многоугольники}.
Иначе говоря, мы допускаем многогранники, вырождающиеся в плоские многоугольники;
под их поверхностью понимается две копии многоугольника, склеенные по границе.

Везде далее мы предполагаем, что многогранник может вырождаться в плоский многоугольник.

\begin{thm}{Теорема Александрова}
\begin{enumerate}[I.]
\item\label{thm:exist} Многогранная метрика на сфере изометрична поверхности выпуклого многогранника тогда и только тогда, когда сумма углов при любой её вершине не превосходит $2\cdot\pi$.

\item\label{thm:unique}
Более того, с точностью до конгруэнтности выпуклый многогранник определяется внутренней метрикой на своей поверхности.
\end{enumerate}

\end{thm}

У Александра Даниловича много замечательных теорем, но с нашей точки зрения эта теорема особенно замечательна.
При этом её доказательство вполне элементарно, его можно объяснить любому, кто знаком с азами топологии.

Эта теорема имеет огромное число приложений.
В частности, она используется в доказательстве её обобщения \cite{alexandrov-1948},
дающего полное описание внутренних метрик на сфере, изометричных замкнутым выпуклым поверхностям  евклидова пространства (не обязательно многогранным).
Последняя теорема является фундаментальным результатом в разделе современной математики --- так называемой \emph{Александровской геометрии}.

Первая часть теоремы является основной, она называется \emph{теоремой существования}.
Вторая часть называется \emph{теоремой единственности}; это лишь небольшая вариация теоремы Коши о многогранниках.
(У Александрова имеется ещё и другая теорема единственности, обобщающая теорему Минковского о многогранниках.)

\begin{wrapfigure}{r}{30mm}
\vskip-4mm
\centering
\includegraphics{mppics/pic-10}
\vskip-0mm
\end{wrapfigure}

Согласно теореме выпуклый многогранник полностью определяется внутренней метрикой своей поверхности.
Значит, зная метрику на поверхности, можно в принципе узнать, где пройдут рёбра многогранника.
Однако на деле это сделать непросто.
Например, поверхность, склеенная из прямоугольника, как показано на картинке, определяет тетраэдр.
При этом некоторые линии склеек оказываются внутри граней тетраэдра и часть рёбер тетраэдра (штрихованные линии) вовсе не идут вдоль сторон прямоугольника.

Теорема доказана А.~Д. Александровым в 1941 году \cite{alexandrov-1941};
мы приведём набросок его доказательства. 
Полное доказательство лучше всего написано в книжке самого Александра Даниловича~\cite{alexandrov}.
Ещё одно доказательство, основанное на деформации трёхмерного многогранного пространства, было предложено Ю. А. Волковым в его кандидатской диссертации \cite{volkov}.

\section{Пространства многогранников и метрик}

\paragraph{Пространство многогранников.}
Обозначим через $\Phi$ пространство всех выпуклых многогранников в евклидовом пространстве, включая многогранники, вырождающиеся в плоский многоугольник.
Многогранники в $\Phi$ будут рассматриваться с точностью до движения пространства, а всё пространство $\Phi$ будет рассматриваться с естественной топологией (достаточно интуитивного понимания близости многогранников друг другу).

Через $\Phi_n$ будем обозначать многогранники в $\Phi$ с ровно $n$ вершинами.
Поскольку любой многогранник в $\Phi$ имеет по крайней мере три вершины,
пространство $\Phi$ разбивается на счётное число подмножеств $\Phi_3,\Phi_4,\dots$

\paragraph{Пространство многогранных метрик.}
Пространство всех многогранных метрик на сфере с суммой углов вокруг каждой точки не более $2\cdot\pi$ будет обозначаться $\Psi$.
При этом метрики в $\Psi$ будут рассматриваться с точностью до изометрии, а всё пространство будет рассматриваться с естественной топологией (опять же, достаточно интуитивного понимания близости двух метрик).

Точка на сфере, вокруг которой сумма углов строго меньше $2\cdot\pi$, будет называться \emph{существенной вершиной}.
Подмножество в $\Psi$, состоящее из метрик с ровно $n$ существенными вершинами, будет обозначаться $\Psi_n$.
Нетрудно видеть, что любая метрика в $\Psi$ имеет по крайней мере три существенные вершины,
так что $\Psi$ разбивается на счётное число подмножеств $\Psi_3,\Psi_4,\dots$

\paragraph{От многогранника к его поверхности.}
Напомним, что поверхность выпуклого многогранника является сферой с многогранной метрикой, такой, что сумма углов вокруг каждой точки не более $2\cdot\pi$.
Таким образом, переход от многогранника к его поверхности определяет отображение 
\[\iota\:\Phi\to \Psi.\]

Заметим также, что число вершин многогранника равно числу существенных вершин его поверхности.
Иначе говоря, $\iota(\Phi_n)\subset \Psi_n$ для любого $n\ge 3$.

\section{О доказательстве}

Используя обозначения, введённые в предыдущей секции, можно дать следующую более точную формулировку теоремы Александрова:

\begin{thm}{Переформулировка теоремы}
Для любого целого $n\ge 3$
отображение $\iota$ является биекцией из $\Phi_n$ в $\Psi_n$.
\end{thm}

Приведённый ниже набросок основан на оригинальном доказательстве А. Д. Александрова.
Главное в нём --- построение однопараметрического семейства многогранников, которое начинается с произвольного многогранника и заканчивается многогранником с поверхностью, изометричной данной.  
Подобный ход рассуждений называется \emph{методом непрерывности};
он часто используется в теории дифференциальных уравнений.

\medskip

Две части первой формулировки теоремы доказываются отдельно.

\parit{Часть \ref{thm:unique}.} Докажем, что отображение $\iota\:\Phi_n\to\Psi_n$ инъективно; иначе говоря, выпуклый многогранник определяется внутренней метрикой на своей поверхности с точностью до движения пространства.

Последнее утверждение является уточнением теоремы Коши о многогранниках. 
Его доказательство практически повторяет доказательство Коши.

Теорема Коши утверждает, что грани многогранника вместе с правилом склейки полностью определяют выпуклый многогранник;
её доказательство приводится во многих классических популярных текстах \cite{aigner-zigler,dolbilin,tabacnikov-fuks}.

\medskip

\parit{Часть \ref{thm:exist}.} Докажем, что отображение $\iota\:\Phi_n\to\Psi_n$ сюръективно.
Это часть доказательства разбивается на следующие леммы.

\begin{thm}{Лемма}
Пространство $\Psi_n$ связно для любого целого $n\ge 3$.
\end{thm}

Доказательство этой леммы несложное, но требует изобретательности; 
его можно провести явным построением непрерывного однопараметрического семейства  метрик в $\Psi_n$, соединяющего две данные метрики.
Такое семейство можно получить путём последовательного применения следующего построения и ему обратного.

Пусть $M$ --- сфера с метрикой из $\Psi_n$, а $v$ и $w$ --- две существенные вершины в $M$.
Разрежем $M$ вдоль кратчайшей линии от $v$ до $w$. 
Заметим, что кратчайшая не проходит через другие существенные вершины в $M$.
Далее заметим, что существует трёхпараметрическое семейство заплаток, которыми можно заклеить разрез так, чтобы полученная метрика осталась в $\Psi_n$; 
в частности, полученная метрика всё ещё будет иметь ровно $n$ существенных вершин (при этом вершины $v$ и $w$ могут перестать быть существенными).

\begin{thm}{Лемма}
Отображение $\iota\:\Phi_n\to\Psi_n$ открыто, 
то есть образ открытого множества в $\Phi_n$ открыт в $\Psi_n$.

В частности, для любого $n\ge 3$, множество $\iota(\Phi_n)$ открыто в~$\Psi_n$.
\end{thm}

Это утверждение очень близко к так называемой \emph{теореме об инвариантности области}, утверждающей, что непрерывное инъективное отображение между многообразиями равной размерности открыто.

Согласно части \ref{thm:exist}, отображение $\iota\:\Phi_n\to\Psi_n$ инъективно.
Доказательство теоремы об инвариантности области теоремы проходит в нашем случае, так как оба пространства $\Phi_n$ и $\Psi_n$ имеют размерность $3\cdot n-6$ и оба похожи на многообразия, хотя, формально говоря, таковыми не являются.
На более техническом языке можно сказать, что оба пространства $\Phi_n$ и $\Psi_n$ имеют естественную структуру $(3\cdot n-6)$-мерного \emph{орбиобразия} и отображение $\iota$ уважает эту структуру.

Мы только убедимся в том, что размерность обоих пространств $\Phi_n$ и $\Psi_n$ равна $3\cdot n-6$.

Выберем многогранник $P$ в $\Phi_n$.
Заметим, что $P$ однозначно определяется $3\cdot n$ координатами своих $n$ вершин.
При этом мы можем считать, что первая вершина совпадает с началом координат, у второй две координаты нулевые, а у третьей одна координата равна нулю; таким образом, для описания многогранников, близких к $P$, достаточно $(3\cdot n-6)$-параметров.
При этом, если многогранник не имеет симметрий, то это описание можно сделать взаимно-однозначным;
в этом случае окрестность точки $P$ в $\Phi_n$ является $(3\cdot n-6)$-мерным многообразием.
Если группа симметрий $P$ нетривиальна,
то данное описание не является взаимно-однозначным, однако это не влияет на размерность $\Phi_n$.

Случай многогранных метрик аналогичен.
Необходимо построить разбиение сферы на плоские треугольники используя только существенные вершины.
По формуле Эйлера, получаем, что число рёбер в данном разбиении равно $3\cdot n-6$.
При этом длины рёбер однозначно определяют метрику и небольшое изменение длин рёбер даёт метрику с теми же свойствами.

\begin{thm}{Лемма}
Отображение $\iota\:\Phi_n\to\Psi_n$ замкнуто, 
то есть образ замкнутого множества в $\Phi_n$ замкнут в $\Psi_n$.

В частности, для любого $n\ge 3$ множество $\iota(\Phi_n)$ замкнуто в~$\Psi_n$.
\end{thm}

Выберем замкнутое множество $Z$ в $\Phi_n$.
Обозначим через $\bar Z$ замыкание множества $Z$ в $\Phi$; заметим, что $Z=\Phi_n\cap \bar Z$.
Предположим, что $P_1,P_2,\dots\in Z$ --- последовательность многогранников, сходящихся к многограннику $P_\infty\in\bar Z$.
Заметим, что $\iota(P_n)$ сходится к $\iota(P_\infty)$  в $\Psi$.
В частности, $\iota(\bar Z)$ замкнуто в~$\Psi$.

Поскольку $\iota(\Phi_n)\subset \Psi_n$ для любого $n\ge 3$ получаем, что  $\iota (Z)=\iota(\bar Z)\cap \Psi_n$, чем и завершается доказательство леммы. 

\medskip

Итак, множество $\iota(\Phi_n)$ --- непусто, замкнуто и открыто в $\Psi_n$, при этом $\Psi_n$ связно для любого $n\ge 3$.
Значит, $\iota(\Phi_n)=\Psi_n$; то есть, $\iota\:\Phi_n\z\to\Psi_n$ сюръективно.
\qeds

\parbf{Благодарности.} Мы благодарим С. Б. Александер, Ю. Д. Бураго и Ж. Цукахара за помощь. 
Статья написана при частичной поддержке РФФИ грант 20-01-00070 и ННФ грант DMS-2005279.

\sloppy
\printbibliography[heading=bibintoc]
\fussy